\documentclass[12pt]{article}
\usepackage{latexsym,epsfig,amssymb,amsmath,amsthm,a4wide}
\usepackage{color} 

\theoremstyle{plain}
\newtheorem*{theorem*}{Theorem}
\newtheorem{thm}{Theorem}
\newtheorem{lem}{Lemma}
\newtheorem{cor}{Corollary}
\newtheorem{prop}{Proposition}

\newcommand{\ebox}{\hfill $\Box$\vspace{2ex}}
\newcommand{\pr}{{\bf Proof.}\ }
\newcommand{\bt}{\begin{thm}}
\newcommand{\et}{\end{thm}}
\newcommand{\bl}{\begin{lem}}
\newcommand{\el}{\end{lem}}
\newcommand{\bp}{\begin{prop}}
\newcommand{\ep}{\end{prop}}
\newcommand{\bc}{\begin{cor}}
\newcommand{\ec}{\end{cor}}
\newcommand{\be}{\begin{eqnarray}}
\newcommand{\ee}{\end{eqnarray}}
\newcommand{\bi}{\begin{itemize}}
\newcommand{\ei}{\end{itemize}}
\newcommand{\beq}{\begin{equation}}
\newcommand{\eeq}{\end{equation}}
\newcommand{\noi}{\noindent}
\newcommand{\Aut}{{\rm Aut}}
\newcommand{\GF}{{\rm GF}}
\newcommand{\PSL}{{\rm PSL}}
\newcommand{\PGL}{{\rm PGL}}
\newcommand{\lcm}{{\rm lcm}}

\newcommand{\ml}{l\kern-0.035cm\char39\kern-0.03cm}

\begin{document}

\title{\vspace{-1.9cm}
Non-orientable regular hypermaps \\ of arbitrary hyperbolic type}

\author{}
\date{}
\maketitle

\begin{center}
\vspace{-1.5cm}

{\large Gareth A. Jones} \\
{\small University of Southampton, Southampton, U.K.}\\

\vspace{4mm}
{\large Martin Ma\v caj} \\
{\small Comenius University, Bratislava, Slovakia}\\

\vspace{4mm}
{\large Jozef \v Sir\'a\v n} \\
{\small Slovak University of Technology, Bratislava, Slovakia, and \\
Open University, Milton Keynes, U.K.}\\

\end{center}
\vskip 7mm

\begin{abstract}

\noi One of the consequences of residual finiteness of triangle groups is that for any given hyperbolic triple  $(\ell,m,n)$ there exist infinitely many regular hypermaps of type $(\ell,m,n)$ on compact orientable surfaces. The same conclusion also follows from a classification of those finite quotients of hyperbolic triangle groups that are isomorphic to linear fractional groups over finite fields. A non-orientable analogue of this, that is, existence of regular hypermaps of a given hyperbolic type on {\em non-orientable} compact surfaces, appears to have been proved only for {\em maps}, which arise when one of the parameters $\ell,m,n$ is equal to $2$.

In this paper we establish a non-orientable version of the above statement in full generality by proving the following much stronger assertion: for every hyperbolic triple $(\ell,m,n)$ there exists an infinite set of primes $p$ of positive Dirichlet density, such that (i) there exists a regular hypermap $\mathcal{H}$ of type $(\ell,m,n)$ on a compact non-orientable surface such that the automorphism group of $\mathcal{H}$ is isomorphic to $\PSL(2,p)$, and, moreover, (ii) the carrier compact surface of {\em every} regular hypermap of type $(\ell,m,n)$ with rotation group isomorphic to $\PSL(2,p)$ is necessarily non-orientable.
\end{abstract}


\medskip


\vskip 3mm
\section{Introduction} \label{sec:intro}

Hypermaps are a generalisation of maps, that is, cellular embeddings of graphs on surfaces, and regular hypermaps are a generalisation of regular maps, in which the map automorphism group is transitive (and hence regular) on mutually incident vertex-edge-face triples. In James' representation \cite{Jam}, a regular hypermap of type $(\ell,m,n)$ is modelled by an embedding of a connected 3-valent graph on a surface, with faces properly 3-coloured and with each face in a monochromatic class bounded, respectively, by $2\ell$-gons, $2m$-gons and $2n$-gons, such that the colour-preserving automorphism group of the embedding is transitive (and hence regular) on vertices of the graph; details will be given in Section~\ref{sec:hyper}.
\smallskip

Every combinatorial structure admitting a regular structure-preserving action of a group on building blocks of the structure may be identified with the group itself; in a similar way a regular hypermap $\mathcal{H}$ of type $(\ell,m,n)$ may also be identified with its automorphism group $\Aut(\mathcal{H})$. This can be done in a particularly appealing way by choosing  three involutory automorphisms $x$, $y$, $z$, acting as reflections in the axes of three edges incident with a fixed vertex in the James' representation of $\mathcal{H}$, so that $\langle x,y,z\rangle = \Aut(\mathcal{H})$. The choice can be made in such a way that the products $r=yz$, $s=zx$ and $t=xy$ act, respectively, as $\ell$-fold, $m$-fold and $n$-fold rotations of the three faces incident with the fixed vertex. The group $\langle r,s,t\rangle$ is also known as the {\em rotation group} of $\mathcal{H}$, and it coincides with the {\em even} subgroup $\Aut_e(\mathcal{H})$ of $\Aut(\mathcal{H})$ consisting of words of even length in the generators $x,y,z$.
\smallskip

It then turns out that the automorphism group $\Aut(\mathcal{H})$ of such a regular hypermap $\mathcal{H}$ of type $(\ell,m,n)$ may be identified with a group of the form $H\langle u\rangle$, where $H = \langle r,s,t\rangle  \ \cong  \Aut_e(\mathcal{H}) $ is a group with presentation
\be\label{eq:pres} \langle r,s,t\ |\ r^{\ell}=s^{m}=t^{n}= rst = \ldots = 1 \rangle \ee
and $u$ is an element of the set $\{x,y,z\}$ inverting two of the three generators of $H$. If the carrier surface of the James' representation of $\mathcal{H}$ is orientable, then $H$ has index $2$ in $H\langle u\rangle$ and forms the group of orientation-preserving automorphisms of $\mathcal{H}$ while the coset $Hu$ comprises the orientation-reversing automorphisms. In the case when the carrier surface of $\mathcal{H}$ is non-orientable, then there is no such distinction and $H=Hu=H\langle u\rangle$. Equivalently, the carrier surface of the hypermap $\mathcal{H}$ is orientable or not, depending on whether the index of $\Aut_e(\mathcal{H})$ in $\Aut(\mathcal{H})$ is $2$ or $1$, and in the orientable case $\Aut_e(\mathcal{H})$ is the group of orientation preserving automorphisms of $\mathcal{H}$. Still another way to look at the situation is to say that $\mathcal{H}$ is {\em inner regular} if $u\in H$, and then inner regularity of $\mathcal{H}$ is equivalent to its non-orientability, cf.~\cite{GAJ}.
\smallskip

Regular hypermaps of type $(\ell,m,n)$ are thus closely related to the well known and widely studied $(\ell,m,n)$-triangle groups $\Delta=\Delta(\ell,m,n)$ with presentation
\be\label{eq:Delta} \Delta = \langle \rho,\sigma,\tau\ | \ \rho^\ell=\sigma^m=\tau^n= \rho \sigma \tau = 1\rangle \ee
and to their extended versions $\Delta\langle \omega\rangle$ obtained by adjoining an involution $\omega$ inverting two of the three generators in $\{\rho,\sigma,\tau\}$. In somewhat more detail, for $\Delta=\Delta (\ell,m,n)$, regular hypermaps of type $(\ell,m,n)$ correspond to quotients of extended triangle groups $\Delta\langle\omega \rangle$, or, equivalently, to their torsion-free normal subgroups, or, still equivalently, to their smooth epimorphic images (those preserving orders of $\rho$, $\sigma$, $\tau$ and $\omega$). Moreover, if $\theta:\  \Delta\langle \omega\rangle\to G$ is such a smooth epimorphism, the hypermap corresponding to $G$ is carried by a non-orientable surface if and only if $\theta(\omega)$ is contained in the $\theta$-image of $\Delta$. These correspondences give rise to rich connections between hypermaps, subgroups of triangle groups, Riemann surfaces and Galois groups; see \cite{JS} for a survey.
\smallskip

Finite regular hypermaps of {\em hyperbolic type} $(\ell,m,n)$, i.e., with $1/\ell+1/m+1/n < 1$, are of particular interest, since the (compact) carrier surfaces of their James' representations have negative Euler characteristic. The fact that for every hyperbolic type there exist infinitely many regular hypermaps of that type on compact {\em orientable} surfaces is a consequence of residual finiteness of the corresponding triangle groups. The same conclusion follows also from a study \cite{Sah} of finite quotients of triangle groups isomorphic to linear fractional groups $\PSL(2,q)$ and $\PGL(2,q)$ for prime powers $q$.
\smallskip

In contrast, very little appears to be known about the existence of hypermaps of a given hyperbolic type but carried by a compact {\em non-orientable} surface. This has generated the following fundamental question in this area of research.
\medskip

\noindent {\bf Question.} {\sl Is it true that for every hyperbolic triple there exist infinitely many finite, non-orientable, regular hypermaps of that type? Equivalently, is it true that for every hyperbolic type $(\ell,m,n)$ the triangle group $\Delta(\ell,m,n)= \langle \rho,\sigma,\tau \ |\ \rho^\ell= \sigma^m= \tau^n =\rho\sigma\tau = 1 \rangle$ contains infinitely many distinct torsion-free normal subgroups $N$ of finite index such that the quotient group $\Delta(\ell,m,n)/N$ contains an involution inverting two of the three generators $\rho N$, $\sigma N$, $\tau N$? }
\bigskip

The aim of this paper is to answer this question in the affirmative by proving the following extremely strong version of the expected result in terms of specifying isomorphism types of automorphism groups of the corresponding hypermaps.
\medskip

\bt\label{t:main}
For every hyperbolic triple $\vartheta = (\ell,m,n)$ there exist an infinite set of primes $\mathcal{S} = \mathcal{S}_\vartheta$ of positive Dirichlet density, such that
\bi
\item[{\rm (i)}] for every $p\in \mathcal{S}$ there exists a non-orientable regular hypermap $\mathcal{H}$ of type $\vartheta$, with $\Aut(\mathcal{H}) \cong \PSL(2,p)$, and, moreover,
\item[{\rm (ii)}] for every $p\in \mathcal{S}$, {\rm every} regular hypermap ${\mathcal H}$ of type $\vartheta$  with rotation group $\Aut_e(\mathcal{H}) \cong \PSL(2,p)$ is non-orientable {\rm(}equivalently, inner regular{\rm)}.
\ei
\et
\medskip

In the very special case of regular {\em maps} (that is, graph embeddings with automorphism group regular on vertex-edge-face triples), which arise when one of the parameters $\ell,m,n$ in a hypermap type is equal to $2$, a weaker version of Theorem \ref{t:main} was established in \cite{JoMS}, with $\PSL(2,q)$ or $\PGL(2,q)$ in place of $\PSL(2,p)$ for a suitable infinite set of prime powers $q$ and with no density claim.
\smallskip

Our proof of Theorem \ref{t:main} in Section \ref{sec:proof} is based on the availability of algebraic tools guaranteeing the existence of a suitable involution inverting two of the three generators in finite epimorphic images of triangle groups in linear fractional groups, on `toggling' this tool between a suitable algebraic number field and a finite field, and on deeper results on the Dirichlet density of sets of primes generating ideals of prime norms in number fields. Refinements of the method are outlined in Section \ref{sec:refin} together with related  remarks.
\smallskip

Although this paper is written using the language of hypermaps, the results and arguments can all be restated in the equivalent language of Riemann and Klein surfaces. Any hypermap $\mathcal H$ of hyperbolic type $\vartheta$ is a quotient of the universal hypermap of type $\vartheta$, on the hyperbolic plane $\mathbb H$, by a subgroup $N$ of the extended triangle group $\Delta\langle\omega\rangle$ of that type. The underlying surface of $\mathcal H$ inherits from $\mathbb H$ the structure of a Riemann surface if $N\le\Delta$, or a Klein surface (non-orientable or with boundary) if $N\not\le\Delta$. In particular, the regular hypermaps constructed in Theorem~\ref{t:main} are on compact non-orientable Klein surfaces without boundary, with $\PSL(2,p)$ acting as a group of automorphisms. Their canonical orientable double covers are compact Riemann surfaces, or equivalently complex projective algebraic curves. Being uniformised by subgroups $N\cap\Delta$ of finite index in the triangle group $\Delta$, these curves are defined over algebraic number fields, by Bely\u\i's Theorem (see~\cite{JS} for details).


\section{Hypermaps, regularity, and existence \\ of regular hypermaps of a given type}\label{sec:hyper}

Hypermaps are usually introduced in one of two equivalent ways, depending on whether the emphasis is on topology or algebra. Apparently, the pioneering work on the topic \cite{Cor1} represents a combination of both ways, which  were subsequently developed in a number of papers out of which we highlight the treatment of \cite{CS} which includes a detailed study of regularity. But the two perhaps most illuminating topological ways of handling hypermaps are by means of bipartite maps \cite{Wal} and trivalent face-3-coloured maps \cite{Jam}, respectively referred to as the {\em Walsh model} and the {\em James model} of a hypermap. In our Introduction we presented the latter and indicated a way this model transforms into an algebraic description of regular hypermaps, those with the highest `level of symmetry'.
\smallskip

For completeness we note that a purely algebraic counterpart to the topological definition of a (not necessarily regular) hypermap of type $(\ell,m,n)$ can be given without considering symmetries -- namely, by  introducing a hypermap as an arbitrary {\em transitive permutation group} with a presentation as in \eqref{eq:pres}. For equivalence of such an approach with the James model and for a large number of more details concerning orientability and other features we refer to \cite{Cor1,CS,Jam}, and for connections of hypermaps to Riemann surfaces, Belyi functions and Galois groups we recommend \cite{JS}.
\smallskip

Returning to regular hypermaps, note that no finiteness assumptions have been made in outlining equivalence of topological and algebraic descriptions. In particular, for $\ell,m,n\ge 2$ the triangle group $\Delta(\ell,m,n)$ with presentation \eqref{eq:Delta} is also a group of orientation-preserving automorphisms of a regular hypermap of type $(\ell,m,n)$, called the {\em universal} hypermap of that type and denoted by ${\mathcal U}(\ell,m,n)$. Its supporting surface is the sphere, the Euclidean plane, or the hyperbolic plane, depending on whether $\frac{1}{\ell}+\frac{1}{m} + \frac{1}{n}$ is greater than, equal to, or less than $1$. Topologically and geometrically, in each of the three cases the hypermap ${\mathcal U}(\ell,m,n)$ is a tessellation of the corresponding simply connected surface formed by properly 3-coloured faces bounded by congruent polygons of lengths $2\ell$, $2m$ and $2m$, three meeting at each vertex. Regularity of ${\mathcal U}(\ell,m,n)$ follows from the existence of an additional orientation-reversing involution not contained in $\Delta(\ell,m,n)$ and inverting two of its generators while preserving the remaining one in \eqref{eq:Delta}.
\smallskip

The fact that regular hypermaps of a given type $(\ell,m,n)$ correspond to quotients of the extended $(\ell,m,n)$-triangle group by torsion-free normal subgroups has a topological counterpart, in that every regular hypermap is smoothly covered by a universal hypermap of the same type. In somewhat more detail, for every regular hypermap ${\mathcal H}$ of type $(\ell,m,n)$ on a carrier surface ${\mathcal S}$ (orientable or not), there exists a smooth regular cover $\theta:\ \tilde{\mathcal S}\to {\mathcal S}$, where $\tilde{\mathcal S}$ is the sphere, the Euclidean plane or the hyperbolic plane (depending on the type of ${\mathcal H}$), such that $\theta$ induces a smooth covering ${\mathcal U}(\ell,m,n)\to {\mathcal H}$.
\smallskip

If a finite regular hypermap ${\mathcal H}$ of type $(\ell,m,n)$ is carried by a compact surface, its Euler characteristic $\chi_{\mathcal H}$ is given by the following consequence of Euler's formula applied to the face-3-coloured representation of the hypermap by the embedded graph $\Gamma_{\mathcal H}$ with $|\Aut({\mathcal H)}|$ vertices:
\begin{equation}\label{eq:hyp}
\frac12\,\, |\Aut({\mathcal H)}|\, \left( 1 - \frac{1}{\ell} - \frac{1}{m} - \frac{1}{n} \right) = - \chi_{\mathcal H}
\end{equation}
It follows that the value of $\mu(\ell,m,n):= \frac{1}{\ell} + \frac{1}{m} + \frac{1}{n}$ is again critical. The type $(\ell,m,n)$ is called {\em elliptic}, {\em parabolic} and {\em hyperbolic}, depending on whether $\mu(\ell,m,n)$ is greater than, equal to, or  less than $1$, respectively. This trichotomy corresponds to the universal hypermap ${\mathcal U}(\ell,m,n)$ being supported by
the sphere, the Euclidean plane or the hyperbolic plane, and also to $\chi_{\mathcal H}$ being positive, zero or negative.
\smallskip

As alluded to in the Introduction, a central question in the theory of finite regular hypermaps is their existence of an arbitrary type. Note first that, due to the obvious symmetry in the presentation \eqref{eq:pres}, when considering types $(\ell,m,n)$ one may assume the ordering $2\le \ell\le m\le n$. With this, for elliptic types the possibilities are $(2,2,n)$ and $(2,3,n)$ for $n\in \{3,4,5\}$, and the only three ordered parabolic types care $(2,3,6)$, $(2,4,4)$ and $(3,3,3)$. In all such cases the existence of regular hypermaps of the corresponding type is well known in the orientable case (on a sphere for all elliptic types, examples being $n$-cycles and their duals and Platonic maps, and on a torus there are infinitely many examples of each of the three parabolic types), and in the non-orientable case on a projective plane for types $(2,2,n)$, $(2,3,4)$ and $(2,3,5)$ but with no regular hypermap on a Klein bottle.
\smallskip

The situation becomes much more interesting when it comes to the existence of regular hypermaps of a given {\em hyperbolic type}. In the orientable case, the fact that for every given hyperbolic type there exists an infinite number of finite regular hypermaps of that type (on an infinite family of compact orientable surfaces) is a consequence of the residual finiteness of hyperbolic triangle groups.
\smallskip

Indeed, since every (infinite) extended triangle group $\Delta^* = \Delta(\ell,m,n)\langle \omega \rangle$ of a hyperbolic type can be embedded in a finite-dimensional matrix group over a number field (see e.g.~\cite{Men} for a particularly appealing embedding in ${\rm SL}_3(K)$ for a number field $K$), by Ma{l\kern-0.035cm\char39\kern-0.03cm}cev's theorem \cite{Mal} the group $\Delta^*$ is residually finite. This means that for every finite subset $S\subset \Delta^*$ not containing the identity the group $\Delta^*$ contains a normal subgroup $N$ of finite index such that $N\cap S=\emptyset$. Since it is known that $\Delta^*$ contains only finitely many conjugacy classes of elements of finite order greater than $1$, we let $S_1$ be a set of distinct representatives of these conjugacy classes. By residual finiteness the group $\Delta^*$ contains a torsion-free normal subgroup $N_1$ of finite index and avoiding $S_1$ and hence (referring to the presentation \eqref{eq:Delta}) disjoint from the set
\[ T = \{\rho^i;\  1\le i\le \ell-1\} \cup \{\sigma^j;\  1\le i\le m-1\} \cup \{\tau^i;\  1\le i\le n-1\} \cup \{\omega\} \ .\]
The quotient $\Delta^*/N_1$ now determines a finite regular hypermap ${\mathcal H}_1$ of type $(\ell,m,n)$ on some compact orientable surface. Moreover, for positive integers $l\ge 2$, letting $u_l$ be an arbitrary non-identity element of $N_{l-1}$ and letting $S_l= S_{l-1}\cup\{u_l\}$, by the same token $\Delta^\ast$ contains a torsion-free normal subgroup $N_l$ of finite index and avoiding $S_l$ (and hence distinct from all the previous subgroups $N_{l'}$ for $l'<l$). As $N_l$ avoids $T$, the quotient $\Delta^\ast/N_l$ defines a finite regular hypermap ${\mathcal H}_l$ of type $(\ell,m,n)$ on some compact orientable surface, distinct from all the previously obtained hypermaps ${\mathcal H}_{l'}$ for $l'<l$. It follows that for every hyperbolic type there are infinitely many finite regular hypermaps of that type. (Explicit calculations based on the embedding derived in \cite{Men} and leading to upper bounds on the order of automorphism groups of the resulting hypermaps can be found in \cite{S,MSI}.)
\smallskip

One can be more specific and ask for the existence of regular hypermaps of a given type on compact orientable surfaces with a given (isomorphism type of) automorphism group. An affirmative answer is known for linear fractional groups $\PSL(2,q)$ and $\PGL(2,q)$ for prime powers $q$, thanks to the classification \cite{Sah} of finite quotients of triangle groups $\Delta (\ell,m,n)$ isomorphic to $\PSL(2,q)$ and $\PGL(2,q)$, with an exact enumeration formula in \cite{Adr} for the groups $\PSL(2,q)$. It is striking that one can go here to the extreme and guarantee that $q$ can be chosen among an infinite set of {\em primes}. Namely, and as an example, it follows from \cite{Sah} and a more detailed follow-up in \cite{CPS} that for every hyperbolic triple $(\ell,m.n)$ and for every prime $p\equiv 1$ mod $2\,\lcm(\ell,m,n)$ there exists a regular hypermap of type $(\ell,m,n)$ on a compact orientable surface, with rotation group isomorphic to $\PSL(2,p)$.
\smallskip

If, however, one hopes to mimic this method and construct finite {\em non-orientable} regular hypermaps of any given hyperbolic type, one faces problems. When considering normal torsion-free subgroups $N$ of finite index in hyperbolic triangle groups $\Delta=\Delta(\ell,m,n)$ given by the presentation \eqref{eq:Delta}, one would have to make sure that the quotient group $\Delta/N$ contains an involution inverting two of its three generators $\rho N$, $\sigma N$ and $\tau N$. It is unclear how this extra property could be guaranteed by residual finiteness, and there appear to be no such attempts in the literature. A difficulty also arises when considering non-orientable regular hypermaps on linear fractional groups over finite fields; by \cite{CPS} there is a necessary and sufficient condition for such a group to contain a suitable inverting involution but it is not obvious how one could use the condition for an arbitrary given hyperbolic type.
\smallskip

This outlined state of the art resulted in the fundamental question on the existence of non-orientable regular hypermaps of arbitrary hyperbolic type, and eventually also in our main result, Theorem \ref{t:main}.


\section{Proof of Theorem \ref{t:main}}\label{sec:proof}

A classification of finite quotients of triangle groups $\Delta(\ell,m,n)$ isomorphic to the linear fractional groups $\PSL(2,q)$ and $\PGL(2,q)$ over finite fields was given in \cite{Sah}, and a more detailed analysis of such groups with applications to the theory of regular hypermaps is available from \cite{CPS}. For our purposes we do not need the full generality of the latter; it turns out that it is sufficient to work with the following facts which are extracted from Propositions 3.2, 4.1 and 6.1 of \cite{CPS} for hyperbolic triples $(\ell,m,n)$ and for Galois fields $\GF(p)$, where the primes $p$ are congruent to $1$ mod $2\,\lcm(\ell,m,n)$; we note that for every given hyperbolic triple there is an infinite set of such primes by Dirichlet's theorem on primes in arithmetic progressions.
\smallskip

\bp \label{p:beg} Let $(\ell,m,n)$ be a hyperbolic triple and let $p$ be a prime such that $p\equiv 1$ {\rm mod} $\,2\,\lcm (\ell, m,n)$. Let $\xi$, $\zeta$ and $\eta$ be, respectively, primitive $(2\ell)^{\rm th}$, $(2m)^{\rm th}$ and $(2n)^{\rm th}$ roots of unity in $\GF(p)$, and let $\xi_\ast = \xi+\xi^{-1}$, $\zeta_\ast = \zeta{+} \zeta^{-1}$ and $\eta_\ast = \eta{+}\eta^{-1}$. Further, assume that the quantity
\be\label{eq:D} d = \xi_\ast^2 + \zeta_\ast^2 + \eta_\ast^2 - \xi_\ast \zeta_\ast \eta_\ast -4
\ee
is not equal to zero. For $\delta=-1/\sqrt{-d}$, which lies in $\GF(p^2)\setminus\GF(p)$ if $-d$ is not a square in $\GF(p)$, and for $\theta=(\xi-\xi^{-1})^{-1}$ let $X$, $Y$ and $Z$ be the elements of $\PSL(2,p)$ or $\PSL(2,p^2)$ defined by
\[ X= \delta\theta\left[ \begin{array}{cc} d & d(\zeta_\ast\xi- \eta_\ast)\\ \eta_\ast-\zeta_\ast\xi^{-1} & -d\end{array}\right], \ \
Y =\delta\left[\begin{array}{cc} 0 & \xi d \\ \xi^{-1} & 0\end{array}\right]
\ \ {\rm and} \ \ Z = \delta\left[\begin{array}{cc} 0 & d \\ 1 & 0 \end{array}\right]\ , \]
with the usual convention that a matrix is identified with its negative. Then, letting $R=YZ$, $S=ZX$ and $T=XY$,
\bi
\item[\rm{(a)}] the group $G=\langle R,S,T\rangle$ is isomorphic to $\PSL(2,p)$ and has a presentation as in {\rm \eqref{eq:pres}} with $R$, $S$ and $T$ in place of $r$, $s$ and $t$; and
\item[\rm{(b)}] the supporting surface of the regular map with automorphism group $\langle X,Y,Z\rangle$ and rotation group $G$ is non-orientable if and only if $-d$ is a square in $\GF(p)$.
\ei
\ep

We begin by outlining our strategy for proving Theorem \ref{t:main}. Given a hyperbolic type $(\ell,m,n)$,
the first step is to construct a number field containing complex primitive $(2j)^{\rm th}$ roots of unity $\alpha$, $\beta$ and $\gamma$ for $j=\ell,m,n$, respectively, together with a square root of the negative of a quantity obtained from \eqref{eq:D} by replacing $\xi$, $\zeta$ and $\eta$ with $\alpha$, $\beta$ and $\gamma$. The second step then consists of using known theorems in algebraic number theory on degrees of lifts of prime ideals to conclude that the ring of integers of the number field constructed in the first step contains an infinite number of ideals of prime norm (for a set of primes of positive Dirichlet density), such that the quotient by each such ideal is a finite field of prime order to which Proposition \ref{p:beg} applies.
\smallskip

At this point it may be appropriate to recall a few facts on ideals in rings of algebraic integers in number fields; to do so we will loosely refer to \cite{AW}. Let $K$ be an algebraic number field, that is, an extension of a finite degree, say, $h$, of the field $\mathbb Q$ of rational numbers; it is well known that $K=\mathbb{Q} (\omega)$ for some $\omega\in K$ with minimal polynomial over $\mathbb{Q}$ of degree $h$. The ring of algebraic integers $\mathcal{O} = \mathcal{O}_K$ in $K$ consists of elements of $K$ that are roots of monic polynomials with integer coefficients. The three important properties of ideals of the ring $\mathcal{O}$ we will need are as follows: (a) every non-zero ideal $J\subset \mathcal{O}$ has a finite index $[\mathcal{O}:J]$, also known as the {\em norm} $N(J)$ of $J$, (b) every non-zero ideal in $\mathcal{O}$ has a unique factorisation into prime ideals, and (c) every prime ideal $J\subset \mathcal{O}$ is maximal, so that in $\mathcal{O}$ maximality of ideals is equivalent to their primality. In the latter case the quotient ring $\mathcal{O}/J$ is a finite field and there exists a unique rational prime $p$ such that $N(J)=p^j$ for some $j\in \{1,2,\ldots,h\}$; the situation when $j=1$, that is, when a prime ideal is of prime norm, will be of special importance for us.
\smallskip

The definition of a norm of an {\em element} of a number field $K$ of degree $h$ over $\mathbb{Q}$ is based on the fact that $K$ admits exactly $h$ distinct injective homomorphisms $\sigma_j:\ K\to \mathbb{C}$, $j\in \{1,2, \ldots, h\}$, into the field $\mathbb{C}$ of complex numbers. For every $z\in K$ the {\em norm} $N(z)$ of $z$ is the product $\sigma_1(z) \sigma_2(z)\ldots \sigma_h(z)$; its constituents $\sigma_j(z)$, $1\le j\le h$, are {\em conjugates of $z$ over $K$}. The norm is a multiplicative function, which means that  $N(z_1z_2)= N(z_1)N(z_2)$ for arbitrary $z_1,z_2\in K$. Norms of elements and ideals in $\mathcal{O}= \mathcal{O}_K$ are related as follows: for every non-zero $z\in \mathcal{O}$ the absolute value of $N(z)$ is equal to the norm of the principal ideal $(z)\subset \mathcal{O}$ generated by $z$. In particular, the norm of every non-zero $z\in \mathcal{O}$ is a non-zero integer, with  $|N(z)|=1$ if and only if $z$ is a unit, i.e., an invertible element of $\mathcal{O}$. Another significant property we will use is that if an element $z\in \mathcal{O}$ belongs to an ideal $J$ in $\mathcal{O}$, then $N(J)$ divides $N(z)$.
\smallskip

Equipped with this we are now ready to proceed to the core of the proof of Theorem~\ref{t:main}. Consider an arbitrary but fixed hyperbolic triple $(\ell,m,n)$. For each $j=\ell,m,n$ let $R_j$ be the set of all the  $\varphi(2j)$ primitive complex $(2j)^{\rm th}$ roots of unity, where $\varphi$ is the Euler totient function; each $R_j$ is known to be the set of roots of the $(2j)^{\rm th}$ cyclotomic polynomial of degree $\varphi(2j)$ and hence an algebraic integer over $\mathbb{Q}$. It follows that for every $\alpha\in R_\ell$, $\beta\in R_m$ and $\gamma\in R_n$ the elements $\alpha_\ast =\alpha+\alpha^{-1}$, $\beta_\ast = \beta + \beta^{-1}$ and $\gamma_\ast = \gamma+\gamma^{-1}$ are also algebraic integers over $\mathbb{Q}$, and hence so is the element
\be\label{eq:A}
d_{(\alpha,\beta,\gamma)}= \alpha_\ast^2 + \beta_\ast^2 + \gamma_\ast^2 - \alpha_\ast \beta_\ast \gamma_\ast -4\ ;
\ee
observe that $d_{(\alpha,\beta,\gamma)}$ is a version of the quantity defined by \eqref{eq:D} but with complex roots of unity $\alpha$, $\beta$ and $\gamma$ in place of $\xi$, $\zeta$ and $\eta$. We now show that a choice of complex roots of unity $\alpha$, $\beta$ and $\gamma$ can always be made to give a non-zero value of $d_{(\alpha,\beta, \gamma)}$.

\bl\label{l:d-non0} For every hyperbolic triple $(\ell,m,n)$ there are primitive $(2\ell)^{\rm th}$, $(2m)^{\rm th}$ and $(2n)^{\rm th}$ complex roots of unity $\alpha$, $\beta$ and $\gamma$, respectively, such that  $d_{(\alpha, \beta,\gamma)}\ne 0$.
\el

\pr By the proof of Proposition 2.3 of \cite{CPS} the element $\alpha^2d_{(\alpha, \beta,\gamma)}$ admits a factorisation of the form
\be\label{eq:d-fact} \alpha^2 d_{(\alpha, \beta,\gamma)} = (\alpha-\beta\gamma)(\alpha - \beta^{-1}\gamma) (\alpha - \beta\gamma^{-1}) (\alpha - \beta^{-1}\gamma^{-1}) \ . \ee
Since there are analogous factorisations with $\alpha^2$ replaced by $\beta^2$ and $\gamma^2$, we may without loss of generality assume that $\ell \ge m\ge n$ in our argument.
\smallskip

The Euler function $\varphi$ has the property that $\varphi(2\ell) > 4$ for every $\ell > 6$. This means that for every $\ell>6$ and for any given choice of $\beta\in R_m$ and $\gamma\in R_n$ one can choose an $\alpha \in R_\ell$ distinct from any of the (at most) four values of the form $\beta^{\pm 1}\gamma^{\pm 1}$, which results in a non-zero  value of $d_{(\alpha, \beta,\gamma)}$ by \eqref{eq:d-fact}. If $6\ge \ell \ge m=n$, then by choosing $\beta=\gamma$ the set $\{ \beta^{\pm 1} \gamma^{\pm 1} \}$ has size $2$ if $m=n=5$ and is a singleton for $m=n\in \{3,4,6\}$, so that there is always a choice of $\alpha\in R_\ell$ giving $d_{(\alpha, \beta,\gamma)}\ne 0$. If $6\ge \ell = m > n$, then for any given $\gamma$ one may take $\alpha=\beta$ such that $\alpha^2\ne \gamma^{\pm 1}$ to force $d_{(\alpha, \beta,\gamma)}\ne 0$. For the remaining seven hyperbolic triples for which $6\ge \ell > m > n$, one has $\gamma^2=1$ and there is always a choice of $\alpha$ and $\beta$ such that $\alpha^2 \ne \beta^{\pm 2}$, eventually giving $d_{(\alpha, \beta,\gamma)}\ne 0$.
\ebox

For our fixed hyperbolic triple $\vartheta = (\ell,m,n)$ let $k = 2\,\lcm(\ell,m,n)$, and let $L = L_\vartheta$ be the number field obtained by adjoining a primitive complex $k^{\rm th}$ root of unity to the field $\mathbb{Q}$ of rational numbers. This is known to be the smallest extension of $\mathbb{Q}$ containing {\em all} primitive $(2\ell)^{\rm th}$, $(2m)^{\rm th}$ and $(2n)^{\rm th}$ complex roots of unity (see, for example, \cite[Theorem 4.27]{Nar}), and it can equivalently be described as the extension of $\mathbb{Q}$ determined by the $k^{\rm th}$ cyclotomic polynomial; the degree of the extension is $[L:\mathbb{Q}] = \varphi(k)$. Further, let $D=D_\vartheta$ denote the set of all the non-zero {\em values} of $d_{(\alpha, \beta,\gamma)}$ for $(\alpha,\beta,\gamma)\in R_\ell\times R_m\times R_n$; observe that $D\ne \emptyset$ by Lemma \ref{l:d-non0}. Finally, let $L^*=L^*_\vartheta = L(\sqrt{-D})$ be the number field obtained from $L$ by adjoining the (finitely many) roots of all the quadratic polynomials $x^2+d$ for $d\in D$. The number field $L^*$ will be our primary object of interest; recall that it depends on the triple $\vartheta$ (and as the triple has been fixed, we will omit the subscript $\vartheta$ in most of what follows). Since $L^*$ arises from $L$ by adjoining square roots, the degree $[L^*:L]$ is a power of two and hence $[L^*:\mathbb{Q}] = 2^e \varphi(k)$ for some non-negative integer $e$. Moreover, as the set of polynomials $x^2+d$ for $d\in D$ is invariant under the action of the Galois group of $L$ over $\mathbb{Q}$, the field $L^*$ is a Galois extension of $\mathbb{Q}$ as well.
\smallskip

The final ingredient we need is a consequence of Corollary 4 of Proposition 7.16 of \cite{Nar} on ideals with prime norm in rings of integers in extensions of algebraic number fields. We state the consequence in a restricted form and in terms introduced above, which is sufficient to be applied to the ring $\mathcal{O}_K$ in our case. For the statement, let us remind the reader that for a given a subset $\mathcal{S}$ of the set $\mathcal{P}$ of all primes, the limit
\[ \lim_{s\to 1^+}\left( \sum_{p\in \mathcal{S}}p^{-s}\right)\Big{/}\left( \sum_{p\in \mathcal{P}}p^{-s}\right) \]
(if it exists) is the {\em Dirichlet density} of $\mathcal{S}$. In most applications, including those considered here, this coincides with the {\em natural density} of $\mathcal S$, the limit (if it exists) as $x\to+\infty$ of the proportion of primes $p\le x$ which are in $\mathcal S$.
\smallskip

\bp{\rm \cite[Corollary 4 of Proposition 7.16]{Nar}}\label{p:Dir}
Let $K$ be a number field that is a Galois extension of $\mathbb{Q}$ of degree $h$, and let $\mathcal{O} = \mathcal{O}_K$ be its ring of algebraic integers. Let $\mathcal{P}_K$ be the set of all rational primes $p$ for which the ideal $p\mathcal{O}$ splits in $K$, that is, $p\mathcal{O}$ factors into a product of distinct prime ideals, each of norm $p$. Then the Dirichlet density of the set $\mathcal{P}_K$ is equal to $1/h$.
\ep

For the fixed hyperbolic type $\vartheta = (\ell,m,n)$ let us now return to the extension $L^*=L(\sqrt{-D})$ of the cyclotomic extension $L$ of $\mathbb{Q}$ of degree $\varphi(k)$ for $k=2\,\lcm(\ell,m,n)$, and to the related concepts and notation introduced earlier, including the rings $\mathcal{O}_L$ and $\mathcal{O}_{L^*}$ of algebraic integers of $L$ and $L^*$. Our aim is to apply Proposition \ref{p:Dir} to the number field $K=L^*\,$ and to the related set $\mathcal{P}_K = \mathcal{P}_{L^*}$ of rational primes $p$, as follows. Let $\mathcal{S}$ be the subset of primes $p\in \mathcal{P}_{L^*}$ with the following two properties:

\bi
\item[{\rm (a)}] the prime $p$ is not a divisor of any of the norms $N(1-\omega^i)$, where $\omega \in R_\ell \cup R_m \cup R_n$ and $i$ is a positive divisor of $2j$ but less than $2j$ for $j=\ell$, $m$, $n$, respectively, and
\item[{\rm (b)}] the prime $p$ is not a divisor of any of the norms $N(d)$ for $d=d_{(\alpha,\beta,\gamma)}\in D$.
\ei

\noindent Clearly, the set $\mathcal{P}_{L^*}\setminus \mathcal{S}$ is finite, and so Proposition \ref{p:Dir} implies that $\mathcal{S}$ has a positive Dirichlet density. We claim that for every $p\in \mathcal{S}$, for every prime ideal $J\subset p\,\mathcal{O}_{L^*}$ of norm $p$ (by Proposition \ref{p:Dir}), and for every triple $(\alpha,\beta,\gamma) \in R_\ell\times R_m\times R_n$ such that $d_{(\alpha, \beta, \gamma)}\in D$, the elements $\alpha+J$, $\beta+J$ and $\gamma+J$ of the finite field $\mathcal{O}_{L^*}/J \cong \GF(p)$ have respective multiplicative orders $2\ell$, $2m$ and $2n$, and all the corresponding elements $d_{(\alpha, \beta, \gamma)} + J$ of the same field are non-zero.
\smallskip

Indeed, suppose without loss of generality that, say, for some $J$ as above the order of $\alpha+J$ is a positive divisor $i$ of $2\ell$ such that $i < 2\ell$. This would mean that $\alpha^i+J = 1+J$, which is equivalent to $1-\alpha^i\in J$. But then the norm $N(J)$ would be a divisor of $N(1-\alpha^i)$ and hence $p$ would divide $N(1-\alpha^i)$, contrary to $p\in \mathcal{P}$. An entirely analogous argument shows that  $d_{(\alpha, \beta,\gamma)} + J \ne J$. As an aside, the fact that orders of $\alpha+J$, $\beta+J$ and $\gamma+J$ in $\mathcal{O}_{L^*}/J \cong \GF(p)$ are $2\ell$, $2m$ and $2n$ implies that every $p\in \mathcal{S}$ is necessarily congruent to $1$ mod \,$2\,\lcm(\ell,m,n)$.
\smallskip

It follows that for each $p\in \mathcal{S}$, each ideal $J\subset p\,\mathcal{O}_{L^*}$ with norm $p$, and for an arbitrary triple $(\alpha,\beta,\gamma) \in R_\ell\times R_m\times R_n$ such that $d_{(\alpha,\beta,\gamma)}\in D$, the elements $\alpha+J$, $\beta+J$ and $\gamma+J$ of the field $\mathcal{O}_{L^*}/J \cong \GF(p)$ may be identified, respectively, with the $(2\ell)^{\rm th}$, $(2m)^{\rm th}$ and $(2n)^{\rm th}$ primitive roots $\xi$, $\zeta$ and $\eta$ mod $p$ from Proposition \ref{p:beg}. Moreover, since the element $d=d_{(\alpha,\beta, \gamma)}$ is such that $d=-u^2$ for some $u\in L^*$ (and hence $u \in \mathcal{O}_{L^*}$, see e.g. \cite[5.1\,C]{Rib}), it follows that $-d$ is a square in $\mathcal{O}_{L^*}/J \cong \GF(p)$.
\smallskip

Our arguments in summary imply that for every $p\in \mathcal{S}$, every regular hypermap of type $\vartheta$ arising from the primitive $(2\ell)^{\rm th}$, $(2m)^{\rm th}$ and $(2n)^{\rm th}$ roots mod $p$ as described in the previous paragraph, is non-orientable, and since $D\ne \emptyset$ there is at least one such hypermap for every $p\in \mathcal{S}$.
\smallskip

This establishes Theorem \ref{t:main} in the case when at least one of the entries $\ell,m,n$ is even, as then for each $p \equiv 1$ mod $2\,\lcm(\ell,m,n)$ there is only one family of hypermaps with rotation group $\PSL(2,p)$, namely, the one described in Proposition \ref{p:beg}. If all of $\ell,m,n$ are odd, however, then for each such prime $p$ there is another family of regular maps of a hyperbolic type $(\ell,m,n)$ with rotation  group $\PSL(2,p)$, cf. \cite{Sah,CPS}. This new family is obtained by taking primitive $\ell^{\rm th}$, $m^{\rm th}$ and $n^{\rm th}$ roots of unity mod $p$ instead, and arises in the same way as indicated in Proposition~\ref{p:beg}. (This family owes its existence to the way a triple of type $\vartheta=(\ell, m, n)$ in $\PSL(2,p)$ can, by suitable sign-changes, be lifted to one of type $2\vartheta$ in ${\rm SL}(2,p)$, unless $\ell, m$ and $n$ are all odd, in which case either $\vartheta$ or $2\vartheta$ is possible, but not both; we refer to details to \cite{Sah} and \cite{CPS}.)
\smallskip

In order to complete the proof of Theorem \ref{t:main} for hypermaps of hyperbolic type $(\ell,m,n)$ with all entries odd and with rotation group $\PSL(2,p)$ for primes $p \equiv 1$ mod $2\,\lcm(\ell,m,n)$, one needs to consider an analogue of the extension of $L$ by adjoining square roots of all possible values of $\sqrt{-D}$ for $D$ containing all non-zero values of $d_{(\alpha,\beta,\gamma)}$ corresponding to both the two families of hypermaps. Since this is a straightforward modification of the previous arguments we leave out the details. \hfill $\Box$


\section{An example}\label{sec:example}

As an illustrative example, based on the results of the first author in \cite{GAJ}, let $\vartheta=(7,2,3)$. Choosing this triple (in any order) means that the resulting non-orientable hypermaps and their orientable double covers  attain the upper bounds of $84(g-2)$ and $84(g-1)$ for the number of automorphisms of a non-orientable and an orientable hypermap (or a Klein and a Riemann surface) of genus $g\ge 3$ and $g\ge 2$ respectively. Equivalently, their common rotation group $\PSL(2,p)$ is both an $H^*$-group and a Hurwitz group, attaining these bounds. Choosing the three parameters in the particular order $(7,2,3)$ means that these hypermaps are in fact maps, of type $\{3,7\}$ in Coxeter's terminology, that is, their faces are triangles and their vertices have valency $7$. Non-orientable maps of this type were considered, nearly 50 years ago, in an unpublished PhD thesis by Wendy Hall~\cite{Hall}, using a different approach, starting with the (orientable) Macbeath--Hurwitz maps with orientation-preserving automorphism group $\PSL(2,p)$ for primes $p\equiv\pm 1$ mod 7 (see~\cite{Macb} for these), and then realising the non-orientable maps as their quotient by an orientation-reversing central involution, when such an element exists. Her necessary and sufficient condition for this, given below, is a special case of our condition that $-d$ should be a square. Here we will apply our general method to the construction of non-orientable hypermaps of type $\vartheta$, comparing our results with those of Hall as we go along.
\smallskip

In dealing with type $\vartheta=(7,2,3)$ there is no choice, in any field, for the traces $\zeta_*$ and $\eta_*$ of the generators of order $2$ and $3$: they must be $0$ and $1$ respectively. However, if $\xi$ is a primitive $14$-th root of unity then there are three choices for the trace $\xi_*$ of the generator of order $7$ in the rotation group, namely $\xi_j=\xi^j+\xi^{-j}$ for $j=1, 3$ and $5$. Then equation~(\ref{eq:D}) becomes
\begin{equation}\label{eq:7}
d=d_j=\xi_j^2-3.
\end{equation}
It is straightforward to check that in any field each $d_j$ is non-zero. Then Proposition~\ref{p:beg} gives generators of $\PSL(2,p)$, and shows that the corresponding hypermap $\mathcal H$  is non-orientable if and only if $-d_j=3-\xi_j^2$ is a square in $\GF(p)$. (This is essentially Theorem 2.9 of~\cite{Hall}.)
\smallskip

At this point Hall, working on a much more restricted problem than that considered in this paper, took a more direct path. Noting that the traces $\xi_j$ are the roots of the polynomial $t^3+t^2-2t-1$, she showed that the elements $-d_j$ are the roots of the polynomial $t^3-4t^2+3t+1$, so their product is $-1$. If $p\equiv 1$ mod 4 then $-1$ is a square in $\GF(p)$, so, because the squares form a subgroup of index $2$ in the multiplicative group, either one or all three of the elements $-d_j$ is a square. We therefore obtain at least one non-orientable hypermap with automorphism group $\PSL(2,p)$ for each prime $p\equiv 1$ or $13$ mod $28$, a set of primes of density $2/\varphi(28)=1/6$. (Some, but not all, primes $p\equiv -1$ or $-13$ mod $28$ also yield two such hypermaps.) Using a similar argument, she obtained the same result for the Macbeath--Hurwitz groups $\PSL(2,p^3)$ for primes $p\equiv 5,9, 17$ or $25$ mod~$28$, a set of primes of density $1/3$.
\smallskip

Returning to the more general method used here, let $L$ be the cyclotomic field generated by the complex $k$-th roots of $1$, where $k=2\,\lcm(7,2,3)=84$, a Galois extension of $\mathbb Q$ of degree $\varphi(84)=24$, and let $L^*$ be the field obtained by adjoining the square roots of the three elements $-d_j$ of $L$ (with $d_j$ referring to \eqref{eq:A} rather than \eqref{eq:7}), so that $L^*$ is a Galois extension of $\mathbb Q$ of degree $h=24\cdot 2^e$ for some non-negative integer $e\le 3$. By Proposition~\ref{p:Dir} and the remarks following it, there is a set $\mathcal S$ of primes, of density $1/h>0$, such that for each $p\in{\mathcal S}$ the group $\PSL(2,p)$ is the automorphism group of a non-orientable hypermap $\mathcal H$ of type $\vartheta$, that is, a map of type $\{3,7\}$.
\smallskip

The density $1/h$ obtained here is less than the density $1/6$ we have obtained by using Hall's approach  in~\cite{Hall}. The reason is that while Hall described the density of primes $p$ for which there exists at least one non-orientable regular map of type $(7,2,3)$ with rotation group isomorphic to $\PSL(2,p^e)$, our Theorem 1 deals with primes $p$ with the property that {\em all} regular maps with rotation group isomorphic to $\PSL(2,p)$ are non-orientable.



\section{Remarks}\label{sec:refin}

The statement of Proposition~\ref{p:Dir}, which played an important role in the proof of Theorem~\ref{t:main}, was deliberately tailored to our needs, and we used it to avoid technicalities familiar to specialists in algebraic number theory. In this setting, Proposition~\ref{p:Dir} may be considered to be the weakest in a series of results describing (in our notation) how the Galois group of  $L^*/\mathbb{Q}$ determines asymptotic properties of prime ideals in any Galois extension $E\le F$ of fields satisfying $\mathbb{Q}\le E\le F\le L^*$. These include the Chebotarev Density Theorem \cite{CH}, which relates sets of primes and their densities to conjugacy classes in the Galois group and their relative sizes.
\smallskip

As regards the asymptotic behaviour of regular hypermaps with rotary group isomorphic to $\PSL(2,q)$ and $\PGL(2,q)$, one may take their description from \cite{Sah} or \cite{CPS} and apply the above outline to the fields $E=\mathbb{Q}(\alpha_\ast, \beta_\ast, \gamma_\ast)$ and $F=E(\sqrt{D'})$ for various subsets $D'\subset D_\vartheta$. An example of such an approach is an analysis of the asymptotic distribution of orientable and non-orientable regular hypermaps of type $(n,2,3)$, based on the groups $\PSL(2,q)$,  presented by the first author in \cite{GAJ}; this builds on and extends the earlier work by Hall~\cite{Hall} on the case $n=7$, mentioned in Section~\ref{sec:example}.
\smallskip

\bigskip
\bigskip

\noindent {\bf Acknowledgement.} Research of the second author was supported by the APVV Research Grant 23-0076 and the VEGA Research Grants 1/0727/22 and 1/0437/23. The third author acknowledges 50\% research support from the EU Recovery and Resilience Plan No. R4:09I03-03-V04-00269, and 50\% support from the APVV Research Grant 22-0005.

\end{document}